\documentclass[a4paper,10pt]{article}

\usepackage[utf8]{inputenc}
\usepackage[T1]{fontenc}
\usepackage[frenchb,english]{babel}
\usepackage{textcomp}
\usepackage{amsmath}
\usepackage{amsfonts}
\usepackage{amssymb}
\usepackage{graphicx}
\usepackage{vmargin}
\usepackage[bookmarks=true,bookmarksnumbered=true,colorlinks=false,pdfborder={0 0 0}]{hyperref}
\usepackage{fancyhdr}
\usepackage{lastpage}
\usepackage{cases}
\usepackage{amsthm}
\newtheorem{definition}{Definition}
\newtheorem{proposition}{Proposition}

\setmarginsrb{25mm}{20mm}{25mm}{15mm}{12pt}{5mm}{0pt}{11mm}

\begin{document}

%% Heading %%
%%%%%%%%%%%%%%%%%%%%%%%%%%%%%%%%%%%%%%%%%%%%%%%%%%%%%%%

\pagestyle{fancy}
\chead{}
\lhead{\textit{B. Lenoir, 2013}}
\rhead{\textit{\thepage/\pageref{LastPage}}}
\cfoot{}
\renewcommand{\headrulewidth}{0pt}
\renewcommand{\footrulewidth}{0pt}

%%%%%%%%%%%%%%%%%%%%%%%%%%%%%%%%%%%%%%%%%%%%%%%%%%%%%%%

%% Title %%
%%%%%%%%%%%%%%%%%%%%%%%%%%%%%%%%%%%%%%%%%%%%%%%%%%%%%%%

\title{A general approach of least squares estimation and optimal filtering}
\author{Benjamin Lenoir\textsuperscript{a} \\ \small \textsuperscript{a} \textit{Onera -- The French Aerospace Lab, 29 avenue de la Division Leclerc, F-92322 Ch\^atillon, France} \\ \\ \small Published in \textit{Optimization and Engineering} \\ \small doi: \href{http://dx.doi.org/10.1007/s11081-013-9217-7}{10.1007/s11081-013-9217-7}}
\date{27 May 2013}
\maketitle

%%%%%%%%%%%%%%%%%%%%%%%%%%%%%%%%%%%%%%%%%%%%%%%%%%%%%%%

%% Abstract %%
%%%%%%%%%%%%%%%%%%%%%%%%%%%%%%%%%%%%%%%%%%%%%%%%%%%%%%%

\begin{abstract}
The least squares method allows fitting parameters of a mathematical model from experimental data. This article proposes a general approach of this method. After introducing the method and giving a formal definition, the transitivity of the method as well as numerical considerations are discussed. Then two particular cases are considered: the usual least squares method and the Generalized Least Squares method. In both cases, the estimator and its variance are characterized in the time domain and in the Fourier domain. Finally, the equivalence of the Generalized Least Squares method and the optimal filtering technique using a matched filter is established.

\end{abstract}

%%%%%%%%%%%%%%%%%%%%%%%%%%%%%%%%%%%%%%%%%%%%%%%%%%%%%%%

%% Article %%
%%%%%%%%%%%%%%%%%%%%%%%%%%%%%%%%%%%%%%%%%%%%%%%%%%%%%%%

\section{Introduction}\label{section:introduction}

The least squares method aims at deriving from experimental data, often plagued by measurement noise, the parameters of a mathematical model describing these data. This method was developed independently by Legendre \cite{legendre1820nouvelles} and Gauss \cite{birkes1993alternative}. The mathematical model allows adding information to the data, which are used to fit unknown parameters. To do so, the data processing minimizes the discrepancy between the data and the model. Therefore, the least squares method can be broadly understood as a minimization problem, which depends on the metric used \cite{el2008least,poch2012orthogonal}.

This article proposes a general approach of the least squares method by considering a general norm deriving from a scalar product. After introducing some general definitions and results, the definition and some properties of the least squares method are detailed. Then, two particular norms are considered and the characterization of the least squares estimator and of its variance are given in the Fourier domain. Finally, the equivalence of a particular least squares method and the optimal filtering technique \cite{pinsker1980optimal,geromel1999optimal} is established. This allows bringing together the temporal approaches and the  frequency approaches aiming at finding a signal in noisy data. In particular, the usual one-dimensional optimal filtering technique used in signal processing is generalized to any dimension.

\section{Framework}

\subsection{Measurement process and parameterization}

Let consider a deterministic time-varying signal $s:\mathbb{R}\rightarrow\mathbb{C}$. The device used to measure this signal makes absolute measurements, i.e. the device introduces no bias, but noise is added in the measurement process. As a result the measurement reads $m(t) = s(t) + e(t)$, where $e$ is a Gaussian stochastic process with a null mean. This process is supposed to be stationary and its Power Spectrum Density, called $S$, as well as its correlation function, called $R$, are supposed to be known\footnote{As stated by the Wiener-Khintchine theorem, $S$ is the Fourier transform of $R$ \cite{lampard1954generalization}.}.

It will be assumed in this article that $N\in\mathbb{N}^*$ measurements are made with a constant time step called $\delta t \in\mathbb{R}_+^*$. The signal is filtered before digitization by a low pass filter with a cut-off frequency equal to $1/(2 \delta t)$ so as to avoid aliasing. In the rest of this article, the notation $\mathbf{y}$ will be a vector of $\mathcal{M}_{N,1}(\mathbb{C})$, whose components are the value of $y$ at the sampling time $k \times \delta t$. It is assumed that the signal $s$ depends linearly of $p$ parameters, called $x_l$ with $l\in\text{\textlbrackdbl}1;p\text{\textrbrackdbl}$, such that
\begin{equation}
 \forall k \in\text{\textlbrackdbl}1;N\text{\textrbrackdbl}, \ \mathbf{s}_k = \sum_{l=1}^p x_l f_l(k \times \delta t),
\end{equation}
where $f_l$ are $p$ functions defined on $\mathbb{R}$. This equation can be written in the matrix form $\mathbf{s} = J\mathbf{x}$, with $J\in\mathcal{M}_{N,p}(\mathbb{C})$ and $J_{kl} = f_l(k \times \delta t) $. Given that measurements are plagued by noise, this leads to
\begin{equation}
 \mathbf{m} = J\mathbf{x} + \mathbf{e}.
\end{equation}
The vector $\mathbf{e}$ is a random vector whose covariance matrix, $\Omega$, is defined by $\Omega_{ij} = R((i-j)\delta t)$. The goal of the least squares method or any identification procedure is to obtain $\mathbf{x}^*$ which is an unbiased linear estimator of $\mathbf{x}$, i.e. of the form $\mathbf{x}^*=A\mathbf{m}$ with $A\in\mathcal{M}_{p,N}$ and $\mathbb{E}[\mathbf{x}^*] = \mathbf{x}$ ($\mathbb{E}$ is the expectation operator). Before going further, the next section gives useful definitions and results concerning discrete time Fourier analysis \cite{yang2009discrete}, generalized to functions defined on a multi-dimensional space.

\subsection{Definitions and general results}\label{subsection:definitions}

Consider $G: \mathbb{Z} \times \text{\textlbrackdbl}1;p\text{\textrbrackdbl} \rightarrow \mathbb{C}$ which is assumed to verify the following property: $\forall l\in\text{\textlbrackdbl}1;p\text{\textrbrackdbl}$, $\sum_{k} |G_{kl}| <\infty$. The Discrete Time Fourier Transform (DTFT) of $G$, called $\mathcal{F}_{\delta t}\{G\}$, is an application defined by
\begin{equation}
   \forall f\in\mathbb{R}, \ \forall l \in \text{\textlbrackdbl}1;p\text{\textrbrackdbl}, \ \mathcal{F}_{\delta t}\{G\}_l(f) = \delta t \sum_{k\in\mathbb{Z}} G_{kl} e^{-i 2 \pi k f \delta t}.
\end{equation}
The inverse DTFT of an application $h:\mathbb{R} \times \text{\textlbrackdbl}1;p\text{\textrbrackdbl} \rightarrow \mathbb{C}$ is an application, called $\mathcal{F}_{\delta t}^{-1}\{h\}$, defined by
\begin{equation}
  \forall k \in \mathbb{Z}, \ \forall l \in \text{\textlbrackdbl}1;p\text{\textrbrackdbl}, \ \mathcal{F}_{\delta t}^{-1}\{h\}_{kl} = \int_{-\frac{1}{2\delta t}}^{-\frac{1}{2 \delta t}} h_l(f) e^{2\pi i k f \delta t} df.
\end{equation}
By extension, the DTFT can be defined for a matrix $M \in \mathcal{M}_{N,p}(\mathbb{R})$. Introducing $\breve{M}: \mathbb{Z} \times \text{\textlbrackdbl}1;p\text{\textrbrackdbl} \rightarrow \mathbb{C}$ defined by $\breve{M}_{kl} = M_{kl}$ if $k\in\text{\textlbrackdbl}1;N\text{\textrbrackdbl}$ and $\breve{M}_{kl} = 0$ otherwise, the DTFT of $M$ is defined by $\mathcal{F}_{\delta t}\{M\} = \mathcal{F}_{\delta t}\{\breve{M}\}$.

Given two applications $A$ and $B$ belonging to $\mathbb{Z} \times \text{\textlbrackdbl}1;p\text{\textrbrackdbl} \rightarrow \mathbb{C}$, the quantity $\left<A|B\right> \in \mathcal{M}_{p p}(\mathbb{C})$ is called the matrix scalar product of $A$ and $B$ and is defined by
\begin{equation}
 \forall (k;l) \in {\text{\textlbrackdbl}1;p\text{\textrbrackdbl}}^2, \ \left<A|B\right>_{kl} = \sum_{j\in\mathbb{Z}} \overline{A}_{jk} B_{jl}.
\end{equation}
The Parseval theorem leads to
\begin{equation}
 \forall (k;l) \in {\text{\textlbrackdbl}1;p\text{\textrbrackdbl}}^2, \ \left<A|B\right>_{kl} = \frac{1}{\delta t} \int_{-\frac{1}{2\delta t}}^{-\frac{1}{2\delta t}} \overline{\mathcal{F}_{\delta t}\{A\}_k(f)} \mathcal{F}_{\delta t}\{B\}_l(f) df.
 \label{eq:Parseval}
\end{equation}

Finally, let $Q$ be an application $\mathbb{Z} \rightarrow \mathbb{C}$. The application $L_Q$, called generalized convolution, is defined by
\begin{equation}
 \forall X:\mathbb{Z}\times\text{\textlbrackdbl}1;p\text{\textrbrackdbl} \rightarrow \mathbb{C}, \ \forall (k;l) \in \mathbb{Z} \times \text{\textlbrackdbl}1;p\text{\textrbrackdbl}, \ L_Q(X)_{kl} = \sum_{i\in\mathbb{Z}} Q_{k-i} X_{il},
\end{equation}
and the notation $L_Q(X) = Q \ast X$ is used. The inverse application of $L_Q$ defined by $L_Q \circ L_Q^{-1} = L_Q^{-1} \circ L_Q = \mathrm{Id}$ has the following expression
\begin{equation}
 L_Q^{-1}(Y) = \mathcal{F}_{\delta t}^{-1} \left\{ \frac{\delta t^2}{\mathcal{F}_{\delta t}\{Q\}} \right\} \ast Y.
\end{equation}

\section{General least-squares solution}

The least squares estimator is the quantity $\mathbf{x}^*$ such that the norm of $\mathbf{m}-J\mathbf{x}$ is minimum for $\mathbf{x}=\mathbf{x}^*$. The result of this computation depends on the choice of the norm on the vector space. In this article, it will be assumed that the norm derives from a scalar product. Therefore regression methods based on a norm not derived from a scalar product, such as the lasso method \cite{osborne2000lasso} \cite{tibshirani1996regression}, are excluded from the discussion. The norm is defined for $\mathbf{y}$ and $\mathbf{z}$ in $\mathcal{M}_{N,1}(\mathbb{C})$ by
\begin{equation}
 \left<\mathbf{y}|\mathbf{z}\right>_P = \mathbf{y}'P\mathbf{z}
\end{equation}
with $\mathbf{y}'$ the conjugate transpose of $\mathbf{y}$ and $P$ a positively defined matrix. The norm is a linear application defined by $|\mathbf{x}|_P = \sqrt{\left<\mathbf{x}|\mathbf{x}\right>_P}$. Since ${|\mathbf{m}-J\mathbf{x}|_P}^2$ is convex in the components of $\mathbf{x}$ and the norm is always positive, the minimum is given by the condition $\nabla_\mathbf{x}\left({|\mathbf{m}-J\mathbf{x}|_P}^2\right) = 2(J'PJ\mathbf{x} - J'P\mathbf{m}) = 0$. Under the condition that $J'PJ$ is invertible, i.e. that $J$ has a rank equal to $p$, the estimator is equal to
\begin{equation}
 \mathbf{x}^* = (J'PJ)^{-1} J' P \mathbf{m}.
 \label{eq:def_estimator}
\end{equation}
Since $\mathbb{E}\left[\mathbf{x}^*\right] = (J'PJ)^{-1}J'P\mathbb{E}\left[J\mathbf{x} + \mathbf{e}\right] = \mathbf{x}$, this estimator is a Gaussian random vector without bias and its covariance matrix $V_{\mathbf{x}^*}$ is equal to
\begin{equation}
 V_{\mathbf{x}^*} = (J'PJ)^{-1}J'P' \Omega PJ(J'PJ)^{-1}.
 \label{eq:def_variance}
\end{equation}

\begin{definition}[Formal definition of the least squares method]
 $N\in\mathbb{N^*}$, $p\in\mathbb{N^*}$ with $p\leq N$. $\mathcal{P}_{N}$ is the ensemble of $N \times N$ symmetric and positively defined matrices. $\mathcal{R}_{N,p}$ is the ensemble of $N \times p$ matrices of rank $p$. The least squares method is defined by the following application
 \begin{eqnarray}
  \mathrm{Ls}:(\mathcal{P}_{N},\mathcal{R}_{N,p},\mathcal{M}_{N,1}) & \rightarrow & \mathcal{M}_{p,1} \\
  (P,J,\mathbf{m}) & \rightarrow & (J'PJ)^{-1}J'P\mathbf{m}
 \end{eqnarray}
\end{definition}

The following result, whose demonstration is straightforward, is of major interest in data processing when consecutive operations using the least squares method are made, for example re-sampling \cite{scott1982effect} or iterative reweighted least squares\cite{holland1977robust}. Indeed, it gives conditions on the scalar products and the projection matrices such that the result of the global optimization is independent of the steps.

\begin{proposition}[Transitivity]
 $(N_0,N_1,N_2)\in(\mathbb{N}^*)^3$ with $N_2 \leq N_1 \leq N_0$. $P_0$, $P_1$ et $\tilde{P}_0$ belong respectively to $\mathcal{P}_{N_0}$, $\mathcal{P}_{N_1}$ and $\mathcal{P}_{N_0}$. $J_0$ and $J_1$ belong respectively to $\mathcal{R}_{N_0,N_1}$ and $\mathcal{R}_{N_1,N_2}$.
 \begin{equation}
  \forall \mathbf{m}  \in \mathcal{M}_{N_0,1} \ , \ \mathrm{Ls}(P_1,J_1,\mathrm{Ls}(P_0,J_0,\mathbf{m})) = \mathrm{Ls}(\tilde{P}_0,J_0 J_1,\mathbf{m})
 \end{equation}
 is equivalent to
 \begin{equation}
  (J_1' P_1 J_1)^{-1} J_1' P_1 (J_0' P_0 J_0)^{-1} J_0' P_0 = (J_1' J_0' \tilde{P}_0 J_0 J_1)^{-1} J_1'J_0'\tilde{P}_0
 \label{eq:transitivity}
 \end{equation}
 \label{prop:transitivity}
\end{proposition}

\begin{proposition}
 $N\in\mathbb{N^*}$, $p\in\mathbb{N^*}$ with $p\leq N$. $J$ belongs to $\mathcal{R}_{N,p}$ and $P$ belongs to $\mathcal{P}_{N}$. $\mathbf{x}^*$ and $V_{\mathbf{x}^*}$ are defined by equations \eqref{eq:def_estimator} and \eqref{eq:def_variance} respectively. Then :
 \begin{itemize}
  \item There exist $\tilde{J} \in \mathcal{R}_{N,p}$ and $C \in \mathcal{M}_{p,p}$ such that $\det(C) \neq 0$, $\tilde{J} = JC$ and $\tilde{J}'P\tilde{J} = \mathrm{Id}$.
  \item $\mathbf{x}^* = C \tilde{J}' P \mathbf{m}$ and $V_{\mathbf{x}^*} = C \tilde{J}' P \Omega P \tilde{J} C'$.
 \end{itemize}
\end{proposition}

The first item of this proposition comes from the Gram-Schmidt algorithm. The second item is interesting since least squares algorithms are to be implemented numerically with the possibility of having $\det(J'PJ)$ close to $0$. This proposition offers the possibility to avoid numerically unstable inversion by replacing it with the Gram-Schmidt process, which does not generate numerical instability \cite{bjorck1967solving,ling1986recursive}. This is therefore an alternative to the ridge regression \cite{hoerl1970ridge,hoerl1975ridge,hoerl1976ridge}.

Out of the multiple possible choices for the matrix $P$, there are two usual ways to define the scalar product, which will be detailed in the next two sections. The most straightforward one is $P = \mathrm{Id}$. The other possibility is $P = {\Omega}^{-1}$.

\section{Usual least squares method}

In this section, $P = \mathrm{Id}$. In this case, $\mathbf{x}^* = (J'J)^{-1} J' \mathbf{m}$ and $V_{\mathbf{x}^*} = (J'J)^{-1}J' \Omega J(J'J)^{-1}$. It is important to notice that this least squares algorithm is not transitive, i.e. it does not fulfill the condition given by equation \eqref{eq:transitivity}, in general. The variance of the estimator can be expressed using the power spectrum density of the measurement noise:
\begin{equation}
 V_{\mathbf{x}^*} = (J'J)^{-1}W(J'J)^{-1}
\end{equation}
with $W\in\mathcal{M}_{p,p}(\mathbb{C})$ defined by
\begin{equation}
 W_{kl} = \int_{-\frac{1}{2 \delta t}}^{\frac{1}{2 \delta t}} S(f) \frac{ \overline{\mathcal{F}_{\delta t}\{J\}_k(f)} \mathcal{F}_{\delta t}\{J\}_l(f) }{T^2}\mathrm{d}f
\end{equation}

\begin{proof}
 The goal is to express $W$ using the PSD of the noise:
 \begin{equation}
  W_{kl} = (J' \Omega J)_{kl} = \sum_{m=1}^{N} \overline{J_{mk}} \sum_{n=1}^{N} R((m-n)\delta t) J_{nl}
 \end{equation}
 The sum over $n$ is a convolution and the sum over $m$ is a scalar product. Expressing these two sums with the Fourier transforms of $J$ and $R$ using the Parseval formula (cf. eq. \eqref{eq:Parseval}) leads to the result.
\end{proof}

\section{Generalized Least Squares}

The other approach, called ``Generalized Least Squares'' (GLS) is defined by $P = {\Omega}^{-1}$. In this case, $\mathbf{x}^* = (J'{\Omega}^{-1}J)^{-1} J'{\Omega}^{-1} \mathbf{m}$ and $V_{\mathbf{x}^*} = (J'{\Omega}^{-1}J)^{-1}$. The Aitken theorem \cite{aitken1935least} states that it is the best linear unbiased estimator (BLUE), but it is not always the only one \cite{luati2011equivalence}. This estimator gives the minimum value to $V_{\mathbf{x}^*}$ \cite{cornillon2007regression}, for the order relation $\leq$ on $\mathcal{M}_{k}(\mathbb{C})$ defined by
\begin{equation}
 A \leq B \Leftrightarrow \exists M\in\mathcal{P}_{k}(\mathbb{C})\cup \{\mathrm{0}_k\} , \ A+M=B.
\end{equation}
Contrary to the previous case, this method is transitive (cf. prop. \ref{prop:transitivity}), which is another advantage in addition to the fact that the covariance of the estimator is minimum.

In this context, it is required to compute the inverse of the covariance matrix $\Omega$. By inverting $\Omega$, one loses information on the correlations between events separated by a time longer than $N\times \delta t$. However, this information is available in the correlation function $R$. To build on this observation, the idea is to extend artificially the size of the vectors up to infinity by filling them with zeros. To do so, let introduce $\breve{J}:\mathbb{Z}\times\text{\textlbrackdbl}1;p\text{\textrbrackdbl} \rightarrow \mathbb{C}$ and $\breve{\mathbf{m}} : \mathbb{Z} \rightarrow \mathbb{C}$ defined by
\begin{subequations}
 \begin{numcases}{\forall(i;j) \in \mathbb{Z} \times \text{\textlbrackdbl}1;p\text{\textrbrackdbl}, \ \breve{J}_{ij} =}
  J_{ij} \ \ \mathrm{if} \ i\in\text{\textlbrackdbl}1;N\text{\textrbrackdbl}\\
  0 \ \ \mathrm{otherwise}
 \end{numcases}
\end{subequations}
\begin{subequations}
 \begin{numcases}{\forall i \in \mathbb{Z} , \ \breve{\mathbf{m}}_i =}
  \mathbf{m}_i \ \ \mathrm{if} \ i\in\text{\textlbrackdbl}1;N\text{\textrbrackdbl}\\
  0 \ \ \mathrm{otherwise}
 \end{numcases}
\end{subequations}
The matrix $\Omega$ being symmetrical, the matrix product $\Omega Y$ with $Y$ a matrix with $N$ lines can be replaced by the generalized convolution product. Computing the inverse of the matrix $\Omega$ when $N\rightarrow\infty$ is equivalent to computing the inverse application of $L_R$ (cf. section \ref{subsection:definitions}). The estimator becomes
\begin{equation}
 \mathbf{x}^* = \left< \breve{J} \middle| {L_R}^{-1}(\breve{J}) \right>^{-1} \left< \breve{J} \middle| {L_R}^{-1}(\mathbf{\breve{m}}) \right>,
\end{equation}
and the variance
\begin{equation}
 V_{\mathbf{x}^*} = \left< \breve{J} \middle| {L_R}^{-1}(\breve{J}) \right>^{-1}.
\end{equation}
This variance is smaller than the one computed with finite size matrices because it is computed optimally in a larger space. And more knowledge on the noise is added in the data processing.

In the case $p=1$, i.e. the data depend on one parameter only, the previous equations can be simplified:
\begin{equation}
 \mathbf{x}^* = \left( \int_{-\frac{1}{2\delta t}}^{-\frac{1}{2\delta t}} \frac{|\mathcal{F}_{\delta t}\{J\}(f)|^2}{S(f)} df \right)^{-1}  \int_{-\frac{1}{2\delta t}}^{-\frac{1}{2\delta t}} \frac{\overline{\mathcal{F}_{\delta t}\{J\}(f)}\mathcal{F}_{\delta t}\{\mathbf{m}(f)\}}{S(f)} df,
 \label{eq:GLS_1D}
\end{equation}
and
\begin{equation}
 V_{\mathbf{x}^*}  = \left( \int_{-\frac{1}{2\delta t}}^{-\frac{1}{2\delta t}} \frac{|\mathcal{F}_{\delta t}\{J\}(f)|^2}{S(f)} df \right)^{-1}.
\end{equation}

\subsection{Influence of the knowledge of the noise Power Spectrum Density}

The approach presented above relies heavily on the knowledge of the noise Power Spectrum Density (PSD), as already underlined \cite{hambaba1992robust}. It is therefore necessary to quantify the error made by using a wrong PSD to process the data. Let call $\mathring{S}$ the PSD used to process the data and $S$ the true PSD. It is assumed that there exists $w:\mathbb{R}\rightarrow\mathbb{C}$ such that
\begin{equation}
 \mathring{S}(f) = S(f) + \epsilon w(f)
\end{equation}
with $\epsilon \ll 1$ and $\forall f\in\mathbb{R}$, $|w(f)|\leq|S(f)|$. For simplicity purpose, $p$ is supposed to be equal to $1$, i.e. $J$ is a column matrix. We call $\mathring{\mathbf{x}}^*$ and $\mathring{V}_{\mathbf{x}^*}$ respectively the estimator and the variance of the GLS method used with $\mathring{S}$. The estimator $\mathring{\mathbf{x}}^*$ is still unbiased. Concerning the variance, it is equal to
\begin{equation}
  \frac{\mathring{V}_{\mathbf{x}^*}-V_{\mathbf{x}^*}}{V_{\mathbf{x}^*}} = \epsilon^2
  \begin{bmatrix}\
  \displaystyle{V_{\mathbf{x}^*} \int_{-\frac{1}{2\delta t}}^{-\frac{1}{2\delta t}} \frac{w^2(f)}{S^3(f)} |\mathcal{F}_{\delta t}\{J\}(f)|^2 df} \\
   \displaystyle{- {V_{\mathbf{x}^*}}^2 \left( \int_{-\frac{1}{2\delta t}}^{-\frac{1}{2\delta t}} \frac{w(f)}{S^2(f)} |\mathcal{F}_{\delta t}\{J\}(f)|^2 df \right)^2}
  \end{bmatrix}
  + o(\epsilon^2),
\end{equation}
which shows that the error is of the order of $\epsilon^2$.

\subsection{Optimal filtering}

Finally, this section shows that the Generalized Least Squares method is equivalent to the optimal filtering technique in the limit of ``infinite'' matrices presented above. The goal is to identify a signal $g:\mathbb{R}\rightarrow\mathbb{C}$ in noise. The optimal filtering technique aims at maximizing the following signal-to-noise ratio at time $t_0$ with respect to $h:\mathbb{R}\rightarrow\mathbb{C}$
\begin{equation}
 \frac{S}{N} = \frac{y_s(t_0)}{\sqrt{\mathbb{E}[|y_e(t_0)|^2]}}
\end{equation}
where $y_s$ and $y_e$ are the convolution of $h$ with $s$ and $e$ respectively. The filter $h$ which maximizes the signal-to-noise ratio, is defined by its Fourier Transform \cite[\S 10.1]{papoulis1977signal}:
\begin{equation}
 \mathcal{F}_{\delta t}\{ h \} (f) = K \frac{\overline{\mathcal{F}_{\delta t}\{ g \}}}{S(f)} e^{-i 2 \pi f t_0}
\end{equation}
where $K$ is a coefficient. $h$ is called a matched filter. By setting
\begin{equation}
 K = \left( \int_{-\frac{1}{2\delta t}}^{-\frac{1}{2\delta t}} \frac{|\mathcal{F}_{\delta t}\{g\}|^2}{S(f)} \right)^{-1},
\end{equation}
one finds the same expression than equation \eqref{eq:GLS_1D}, except for $t_0$ which is an offset. This demonstrate the equivalence between the optimal filtering technique and the GLS method when taking ``infinite'' vectors. The GLS method has the advantage to fit several parameters from the data, i.e. to identify several pattern from the data. It also assesses the independence of the fitted parameters by computing their covariance matrix. On the contrary, the usual optimal filtering method can only look for one pattern.

\section{Conclusion}\label{section:conclusion}

This article presented a general approach of least squares estimation. A general definition of this minimization method was introduced and general results were discussed: the conditions required to have the transitivity of the method and a way to avoid numerical instability. Then two particular cases were considered. In both cases, the estimator and its variance were characterized in the Fourier domain. Finally, the equivalence between the Generalized Least Squares method and the optimal filtering technique was established, bridging the gap between signal analysis and optimization.

\section*{Acknowledgments}\label{section:acknowledgments}

The author is grateful to CNES (Centre National d'\'Etudes Spatiales) for its financial support.

%%%%%%%%%%%%%%%%%%%%%%%%%%%%%%%%%%%%%%%%%%%%%%%%%%%%%%%

%% Bibliography %%
%%%%%%%%%%%%%%%%%%%%%%%%%%%%%%%%%%%%%%%%%%%%%%%%%%%%%%%

%%%%%%%%%%%%%%%%%%%%%%%%%%%%%%%%%%%%%%%%%%%%%%%%%%%%%%%


\begin{thebibliography}{10}

\bibitem{aitken1935least}
A.~C. Aitken.
\newblock {On least squares and linear combinations of observations}.
\newblock {\em Proceedings of the Royal Society of Edinburgh}, 55:42--48, 1935.

\bibitem{birkes1993alternative}
D.~Birkes and Y.~Dodge.
\newblock {\em {Alternative methods of regression}}.
\newblock Wiley Series in Probability and Statistics. Wiley, New York, 1993.

\bibitem{bjorck1967solving}
{\AA}.~Bj{\"o}rck.
\newblock {Solving linear least squares problems by Gram-Schmidt
  orthogonalization}.
\newblock {\em BIT Numerical Mathematics}, 7(1):1--21, 1967.

\bibitem{cornillon2007regression}
P.~A. Cornillon and {\'E}.~Matzner-L{\o}ber.
\newblock {\em {R\'egression: Th\'eorie et applications}}.
\newblock Springer Paris, 2007.

\bibitem{el2008least}
M.~I. El-Khaiary.
\newblock {Least-squares regression of adsorption equilibrium data: comparing
  the options}.
\newblock {\em J. Hazard. Mater.}, 158(1):73--87, 2008.

\bibitem{geromel1999optimal}
J.~C. Geromel.
\newblock {Optimal linear filtering under parameter uncertainty}.
\newblock {\em Signal Processing, IEEE Transactions on}, 47(1):168--175, 1999.

\bibitem{hambaba1992robust}
M.~L. Hambaba.
\newblock {The robust generalized least-squares estimator}.
\newblock {\em Signal Process.}, 26(3):359--368, 1992.

\bibitem{hoerl1975ridge}
A.~E. Hoerl, R.~W. Kannard, and K.~F Baldwin.
\newblock {Ridge regression: some simulations}.
\newblock {\em Communications in Statistics}, 4(2):105--123, 1975.

\bibitem{hoerl1970ridge}
A.~E. Hoerl and R.~W. Kennard.
\newblock {Ridge regression: Biased estimation for nonorthogonal problems}.
\newblock {\em Technometrics}, 12(1):55--67, 1970.

\bibitem{hoerl1976ridge}
A.~E. Hoerl and R.~W. Kennard.
\newblock {Ridge regression iterative estimation of the biasing parameter}.
\newblock {\em Communications in Statistics-Theory and Methods}, 5(1):77--88,
  1976.

\bibitem{holland1977robust}
P.~W. Holland and R.~E. Welsch.
\newblock {Robust regression using iteratively reweighted least-squares}.
\newblock {\em Communications in Statistics-Theory and Methods}, 6(9):813--827,
  1977.

\bibitem{lampard1954generalization}
D.~G. Lampard.
\newblock {Generalization of the Wiener-Khintchine Theorem to Nonstationary
  Processes}.
\newblock {\em J. Appl. Phys.}, 25(6):802--803, 1954.

\bibitem{legendre1820nouvelles}
A.~M. Legendre.
\newblock {\em {Nouvelles m\'ethodes pour la d\'etermination des orbites des
  com\`etes}}.
\newblock Chez Firmin Didot, 116 rue de Thionville, Paris, 1820.

\bibitem{ling1986recursive}
F.~Ling, D.~Manolakis, and J.~Proakis.
\newblock {A recursive modified Gram-Schmidt algorithm for least-squares
  estimation}.
\newblock {\em IEEE Trans. Acoust. Speech Signal Process.}, 34(4):829--836,
  1986.

\bibitem{luati2011equivalence}
A.~Luati and T.~Proietti.
\newblock {On the equivalence of the weighted least squares and the generalised
  least squares estimators, with applications to kernel smoothing}.
\newblock {\em Ann. Inst. Stat. Math.}, 63(4):851--871, 2011.

\bibitem{osborne2000lasso}
M.~R. Osborne, B.~Presnell, and B.~A. Turlach.
\newblock On the lasso and its dual.
\newblock {\em Journal of Computational and Graphical Statistics},
  9(2):319--337, 2000.

\bibitem{papoulis1977signal}
A.~Papoulis.
\newblock {\em {Signal analysis}}.
\newblock McGraw-Hill, 1977.

\bibitem{pinsker1980optimal}
M.~S. Pinsker.
\newblock Optimal filtering of square-integrable signals in gaussian noise.
\newblock {\em Probl. Peredachi Inf.}, 16(2):52--68, 1980.

\bibitem{poch2012orthogonal}
J.~Poch and I.~Villaescusa.
\newblock {Orthogonal Distance Regression: A Good Alternative to Least Squares
  for Modeling Sorption Data}.
\newblock {\em J. Chem. Eng. Data}, 57(2):490--499, 2012.

\bibitem{scott1982effect}
A.~J. Scott and D.~Holt.
\newblock {The effect of two-stage sampling on ordinary least squares methods}.
\newblock {\em J. Am. Stat. Assoc.}, 77(380):848--854, 1982.

\bibitem{tibshirani1996regression}
R.~Tibshirani.
\newblock {Regression shrinkage and selection via the lasso}.
\newblock {\em Journal of the Royal Statistical Society. Series B
  (Methodological)}, 58(1):267--288, 1996.

\bibitem{yang2009discrete}
Won~Young Yang.
\newblock {Discrete-Time Fourier Analysis}.
\newblock In {\em Signals and Systems with MATLAB}, pages 129--205. Springer
  Berlin Heidelberg, 2009.

\end{thebibliography}
\end{document}